\documentclass{amsart}
\usepackage{amssymb}
\usepackage{amsmath}
\usepackage{amsfonts}

\newtheorem{theorem}{Theorem}
\theoremstyle{plain}

\newtheorem{corollary}{Corollary}

\newtheorem{remark}{Remark}

\numberwithin{equation}{section}
\input{tcilatex}

\begin{document}
\title[Variation Embedding Theorem]{A Variation Embedding Theorem and
Applications}
\author{Peter Friz, Nicolas Victoir}

\begin{abstract}
Fractional Sobolev spaces, also known as Besov or Slobodetzki spaces, arise
in many areas of analysis, stochastic analysis in particular. We prove an
embedding into certain $q$-variation spaces and discuss a few applications.
First we show $q$-variation regularity of Cameron-Martin paths associated to
fractional Brownian motion and other Volterra processes. This is useful, for
instance, to establish large deviations for enhanced fractional Brownian
motion. Second, the $q$-variation embedding, combined with results of rough
path theory, provides a different route to a regularity result for
stochastic differential equations by Kusuoka. Third, the embedding theorem
works in a non-commutative setting and can be used to establish H\"{o}%
lder/variation regularity of rough paths.
\end{abstract}

\maketitle

\section{Fractional Sobolev Spaces}

For a real valued measurable path $h:\left[ 0,1\right] \rightarrow \mathbb{R}
$ and $\delta \in \left( 0,1\right) $ and $p\in \left( 1,\infty \right) $ we
define the fractional Sobolev (semi-)norm%
\begin{equation*}
\left\vert h\right\vert _{W^{\delta ,p}}=\left( \int \int_{\left[ 0,1\right]
^{2}}\frac{\left\vert h_{t}-h_{s}\right\vert ^{p}}{\left\vert t-s\right\vert
^{1+\delta p}}dsdt\right) ^{1/p}\in \lbrack 0,+\infty ]
\end{equation*}%
For $\delta =1$ and $\,p\in \left( 1,\infty \right) $, writing $\dot{h}$ for
the weak derivative, we set%
\begin{equation*}
\left\vert h\right\vert _{W^{1,p}}=\left( \int_{0}^{1}\left\vert \dot{h}%
_{t}\right\vert ^{p}dt\right) ^{1/p}\in \lbrack 0,+\infty ].
\end{equation*}%
Define $W^{\delta ,p}$ as the set of $h$ for which $\left\vert h\right\vert
_{L^{p}}+\left\vert h\right\vert _{W^{\delta ,p}}<\infty $. They are known
to be Banach-spaces. For $1\geq \delta >1/p>0$ one can assume that $h$ is
continuous; compare with the embedding theorems below. It then makes sense
to consider the closed subspace%
\begin{equation*}
W_{0}^{\delta ,p}=\left\{ h\in W^{\delta ,p}:h\left( 0\right) =0\right\} 
\end{equation*}%
which is Banach under $\left\vert \cdot \right\vert _{W^{\delta ,p}}$. We
finally remark that the space $W^{1,p}$ is precisely the set of absolutely
continuous paths on $\left[ 0,1\right] $ with (a.e. defined) derivative in $%
L^{p}\left[ 0,1\right] $. The space $W_{0}^{1,2}$ is the usual
Cameron-Martin space for Brownian motion. We recall some well-known
continuous resp. compact embeddings\footnote{%
The symbol $\subset \subset $ means compact embedding.} \cite{Ad75},\cite%
{FdeP99},\cite{De04},%
\begin{eqnarray}
p &\in &\left( 1,\infty \right) ,\,\,1\geq \tilde{\delta}>\delta \geq
0\implies W^{\tilde{\delta},p}\subset \subset W^{\delta ,p}\text{,}
\label{CompactEmbeddingWdeltap} \\
1 &<&p\leq q<\infty ,\,\delta \equiv 1-1/p+1/q>0\implies W^{1,p}\subset
W^{\delta ,q}.  \label{W1pEmbedding}
\end{eqnarray}

\section{\protect\bigskip A $q$-Variation Embedding}

\begin{theorem}
Let $p\in \left( 1,\infty \right) $ and $\alpha =1-1/p>0$. Then the
variation of any $h\in W^{1,p}$ is controlled by the control function%
\footnote{%
A continuous, super-additive map $\left( s,t\right) \mapsto \omega \left(
s,t\right) \in \lbrack 0,\infty ),$ defined for $0\leq s\leq t\leq 1$.}%
\begin{equation*}
\omega \left( s,t\right) =\left\vert h\right\vert _{W^{1,p};\left[ s,t\right]
}\left( t-s\right) ^{\alpha },\,\,\,0\leq s\leq t\leq 1
\end{equation*}%
and we have the continuous embeddings%
\begin{equation*}
W^{1,p}\subset C^{\alpha \text{-H\"{o}lder}}\text{ \ and \ \ }W^{1,p}\subset
C^{1\text{-var}}\text{.}
\end{equation*}
\end{theorem}

\begin{proof}
By absolute continuity and H\"{o}lder's inequality with conjugate exponents $%
p$ and $1/\alpha $ 
\begin{eqnarray*}
\left\vert h_{s,t}\right\vert &=&\int_{s}^{t}\left\vert \dot{h}%
_{r}\right\vert dr\leq \left( t-s\right) ^{\alpha }\left(
\int_{s}^{t}\left\vert \dot{h}_{r}\right\vert ^{p}dr\right) ^{1/p} \\
&=&\left\vert h\right\vert _{W^{1,p};\left[ s,t\right] }\left( t-s\right)
^{\alpha }.
\end{eqnarray*}%
We now show that the variation of $h$ is controlled by the control function%
\begin{equation*}
\omega \left( s,t\right) =\left\vert h\right\vert _{W^{1,p};\left[ s,t\right]
}\left( t-s\right) ^{\alpha },\,\,\,t\geq s.
\end{equation*}%
Only super-additivity, $\omega \left( s,t\right) +\omega \left( t,u\right)
\leq \omega \left( s,u\right) $ with $s\leq t\leq u$, is non-trivial. Note $%
p\in \left( 1,\infty \right) $. From H\"{o}lder's inequality with conjugate
exponents $p$ and $p/\left( p-1\right) =1/\alpha $ we obtain%
\begin{eqnarray*}
&&\left\vert h\right\vert _{W^{1,p};\left[ s,t\right] }\left( t-s\right)
^{\alpha }+\left\vert h\right\vert _{W^{1,p};\left[ t,u\right] }\left(
u-t\right) ^{\alpha } \\
&\leq &\left( \left\vert h\right\vert _{W^{1,p};\left[ s,t\right]
}^{p}+\left\vert h\right\vert _{W^{1,p};\left[ t,u\right] }^{p}\right) ^{1/p}%
\left[ \left( t-s\right) ^{\alpha \frac{p}{p-1}}+\left( u-t\right) ^{\alpha 
\frac{p}{p-1}}\right] ^{\left( p-1\right) /p} \\
&=&\left\vert h\right\vert _{W^{1,p};\left[ s,u\right] }\left( u-t\right)
^{\alpha }.
\end{eqnarray*}%
This shows that $\omega $ is super-additive and we conclude that for any $%
0\leq a<b\leq 1,$%
\begin{equation*}
\left\vert h\right\vert _{1\text{-var;}[a,b]}\leq \omega \left( a,b\right)
=\left\vert b-a\right\vert ^{\alpha }\left\vert h\right\vert _{W^{1,p};\left[
a,b\right] }.
\end{equation*}%
In particular, we established $W^{1,p}\subset C^{\alpha \text{-H\"{o}lder}}$%
\ and$\ \ W^{1,p}\subset C^{1\text{-var}}.$
\end{proof}

\begin{theorem}
\label{ThmfroDeltaLTone}Let $0<\delta <1$ and $p\geq 1$ such that%
\begin{equation*}
\alpha =\delta -1/p>0\text{.}
\end{equation*}%
Set $q=1/\delta $. Then the $q$-variation of any $h\in W^{\delta ,p}$ is
controlled by a constant multiple of the control function%
\begin{equation*}
\omega \left( s,t\right) =\left\vert h\right\vert _{W^{\delta ,p};\left[ s,t%
\right] }^{q}\left( t-s\right) ^{\alpha q},\,\,\,0\leq s\leq t\leq 1\text{.}
\end{equation*}%
and we have the continuous embeddings%
\begin{equation*}
W^{\delta ,p}\subset C^{\alpha \text{-H\"{o}lder}}\text{ \ and \ \ }%
W^{\delta ,p}\subset C^{q\text{-var}}\text{.}
\end{equation*}
\end{theorem}

\begin{proof}
We have%
\begin{equation*}
\left\vert h\right\vert _{W^{\delta ,p};\left[ s,t\right] }^{p}\equiv
F_{s,t}=\int \int_{[s,t]^{2}}\frac{\left\vert h_{u,v}\right\vert ^{p}}{%
\left\vert v-u\right\vert ^{1+\delta p}}dudv=\int \int_{[s,t]^{2}}\left( 
\frac{\left\vert h_{u,v}\right\vert }{\left\vert v-u\right\vert ^{1/p+\delta
}}\right) ^{p}dudv.
\end{equation*}%
The Garsia-Rodemich-Rumsey lemma with $\Psi \left( \cdot \right) =\left(
\cdot \right) ^{p}$ and $p\left( \cdot \right) =\left( \cdot \right)
^{1/p+\delta }$ yields%
\begin{eqnarray*}
\left\vert h_{s,t}\right\vert &\leq &C\int_{0}^{t-s}\left( \frac{F_{s,t}}{%
u^{2}}\right) ^{1/p}dp\left( u\right) =C\left\vert h\right\vert _{W^{\delta
,p};\left[ s,t\right] }\int_{0}^{t-s}u^{-2/p}dp\left( u\right) \\
&=&C\left\vert h\right\vert _{W^{\delta ,p};\left[ s,t\right]
}\int_{0}^{t-s}u^{-1/p+\delta -1}du=C\left\vert h\right\vert _{W^{\delta ,p};%
\left[ s,t\right] }\left( t-s\right) ^{\delta -1/p},
\end{eqnarray*}%
using $\alpha \equiv \delta -1/p>0$. We now show that the $q$-variation of $%
h $ is controlled by the control function%
\begin{equation*}
\omega \left( s,t\right) :=\left\vert h\right\vert _{W^{\delta ,p};\left[ s,t%
\right] }^{q}\left( t-s\right) ^{\alpha q},\,\,\,t\geq s.
\end{equation*}%
Only super-additivity, $\omega \left( s,t\right) +\omega \left( t,u\right)
\leq \omega \left( s,u\right) $ with $s\leq t\leq u$, is non-trivial. Note
that $p/q=1/\left( p\alpha +1\right) \in \left( 1,\infty \right) $. From H%
\"{o}lder's inequality with conjugate exponents $p/q$ and $p/\left(
p-q\right) $ we obtain%
\begin{eqnarray*}
&&\left\vert h\right\vert _{W^{\delta ,p};\left[ s,t\right] }^{q}\left(
t-s\right) ^{q\alpha }+\left\vert h\right\vert _{W^{\delta ,p};\left[ t,u%
\right] }^{q}\left( u-t\right) ^{q\alpha } \\
&\leq &\left( \left\vert h\right\vert _{W^{\delta ,p};\left[ s,t\right]
}^{p}+\left\vert h\right\vert _{W^{\delta ,p};\left[ t,u\right] }^{p}\right)
^{q/p}\left[ \left( t-s\right) ^{q\alpha \frac{p}{p-q}}+\left( u-t\right)
^{q\alpha \frac{p}{p-q}}\right] ^{\left( p-q\right) /p}
\end{eqnarray*}%
The first factor is easily estimated 
\begin{equation*}
\left( \left\vert h\right\vert _{W^{\delta ,p};\left[ s,t\right]
}^{p}+\left\vert h\right\vert _{W^{\delta ,p};\left[ t,u\right] }^{p}\right)
^{q/p}\leq \left\vert h\right\vert _{W^{\delta ,p};\left[ s,u\right] }^{q}.
\end{equation*}%
To estimate the second factor note that the exponent of $\left( t-s\right) $
resp. $\left( u-t\right) $ equals one, indeed%
\begin{equation*}
q\alpha \frac{p}{p-q}=1\Longleftrightarrow q=\frac{p}{p\alpha +1}
\end{equation*}%
and the second factor equals%
\begin{equation*}
\left( u-s\right) ^{\left( p-q\right) /p}=\left( u-s\right) ^{q\alpha }.
\end{equation*}%
This shows that $\omega $ is super-additive and we conclude that for any $%
0\leq a<b\leq 1,$%
\begin{equation*}
\left\vert h\right\vert _{q\text{-var;}[a,b]}\leq C\omega \left( a,b\right)
^{1/q}=C\left\vert b-a\right\vert ^{\alpha }\left\vert h\right\vert
_{W^{\delta ,p};\left[ a,b\right] }.
\end{equation*}%
In particular, we have established continuity of the embeddings%
\begin{equation*}
W^{\delta ,p}\subset C^{\alpha \text{-H\"{o}lder}}\ \text{and\ }\ W^{\delta
,p}\subset C^{q\text{-var}}.
\end{equation*}
\end{proof}

The case $p=2$ deserves special attention. The assumptions of Theorem \ref%
{ThmfroDeltaLTone} are then satisfied for any $\delta \in \left(
1/2,1\right) $.

\begin{remark}
In \cite{Ku}, Kusuoka discusses differentiability of SDE solution beyond the
usual Malliavin sense. In particular, he shows the existence of a nice
version of the It\^{o}-map which has derivatives in directions $%
W_{0}^{\delta ,2}\supset W_{0}^{1,2}$ for $\delta \in \left( 1/2,1\right) $.
Since $W_{0}^{\delta ,2}\subset C^{q\text{-var}}$ with $q=1/\delta <2$ this
result is now explained by Lyons' theory of rough paths \cite{LQ97, L98}.
Note that in Lyons' continuity statements the modulus $\omega $ is
preserved. This implies that after perturbating a Brownian path in a $%
W_{0}^{\delta ,2}$-direction the solution maintaines $\alpha $-H\"{o}lder
regularity with $\alpha =\delta -1/2$. (Clearly, this is not true for an
arbitrary perturbation in $C^{q\text{-var}}$!) For what it's worth, we can
extend Gateaux-differentiabiliy to suited $W_{0}^{\delta ,p}$-spaces as long
as $\delta -1/p>0$ and even apply this to rough path differential equations
driven by enhanced fBM. On the other hand, we do not attempt to recover
Kusuoka's full statement (Fr\'{e}chet in both starting point and $%
W_{0}^{\delta ,2}$). This requires a careful formulation of Lyons' universal
limit theorem and will be addressed in a forthcoming monograph.
\end{remark}

\begin{remark}
Integrals of form $\int fdg$ for $f,g\in W^{\delta ,2}$ are discussed in 
\cite{Za01}. Theorem \ref{ThmfroDeltaLTone} reveals them as normal
Young-intergal. Following \cite{LQ02} its continuity properties are
conveniently expressed in terms of the modulus $\omega $. In particular, the
modulus of continuity of $\int fdg$ is immediately controlled by the $%
W^{\delta ,2}$-Sobolev-norms\ of $f$ and $g$ and we can easily extend this
to $W^{\delta ,p}$ provided $\delta -1/p>0$. On the other hand, our approach
does not allow us to control the $W^{\delta ,2}$-norm of the indefinite
integral $\int fdg$.
\end{remark}

\begin{remark}
\label{GausstailForEBM}When $\delta <1$ the notion of $W^{\delta ,p}$ makes
perfect sense for paths with values in a metric space $\left( E,d\right) $.
Theorem \ref{ThmfroDeltaLTone} still holds with the same proof\footnote{%
Simply write $h_{s,t}\equiv d\left( h_{s},h_{t}\right) $ and note that the
Garsia-Rodemich-Rumsey lemma works for $\left( E,d\right) $-valued
continuous functions.}. The case of the free step-$N$ nilpotent group $%
\left( G^{N}\left( \mathbb{R}^{d}\right) ,\otimes \right) $ with
Carnot-Caratheodory norm $\left\Vert \cdot \right\Vert $ and distance $%
d\left( x,y\right) =\left\Vert x^{-1}\otimes y\right\Vert $ is of particular
importance: Theorem \ref{ThmfroDeltaLTone} is a criterion for variation and H%
\"{o}lder regularity of a $G^{N}\left( \mathbb{R}^{d}\right) $-valued path,
a fundamental aspect in Lyons' theory of rough paths, \cite{L98}. To
illustrate the idea we give a simple application to enhanced Brownian motion 
$\mathbf{B}$, see \cite{FV1, FLS}. Then\footnote{%
Note $\left\Vert \mathbf{B}_{s,t}\right\Vert \overset{\mathcal{D}}{=}%
\left\vert t-s\right\vert ^{1/2}\left\Vert \mathbf{B}_{0,1}\right\Vert $.}%
\begin{equation*}
\mathbb{E}\left\Vert \mathbf{B}\right\Vert _{W^{\delta ,p};\left[ 0,1\right]
}^{p}=\int \int_{[0,1]^{2}}\frac{\mathbb{E}\left\Vert \mathbf{B}%
_{s,t}\right\Vert ^{p}}{\left\vert t-s\right\vert ^{1+\delta p}}dsdt=\mathbb{%
E}\left\Vert \mathbf{B}_{0,1}\right\Vert ^{p}\int \int_{[0,1]^{2}}\left\vert
t-s\right\vert ^{p/2-1-\delta p}dsdt.
\end{equation*}%
For every $\alpha <1/2$ and $\delta \in \left( \alpha ,1/2\right) $ there
exists $p_{0}\left( \delta \right) $ such that for all $p\geq p_{0}$ the
double integral is bounded by $1$. Thus for all $p$ large enough,%
\begin{equation*}
\mathbb{E}\left\Vert \mathbf{B}\right\Vert _{W^{\delta ,p};\left[ 0,1\right]
}^{p}\leq \,\mathbb{E}\left\Vert \mathbf{B}_{0,1}\right\Vert ^{p}\text{.}
\end{equation*}%
Is is well-known, \cite{FV1}, that $\left\Vert \mathbf{B}_{0,1}\right\Vert $
has a Gaussian tail and it follows that $\left\Vert \mathbf{B}\right\Vert
_{W^{\delta ,p}}$ has a Gaussian tail, provided $p\geq p_{0}\left( \delta
\right) $. For $p$ large enough we have $\alpha \leq \delta -1/p$ and we
conclude that $\left\Vert \mathbf{B}\right\Vert _{\alpha \text{-H\"{o}lder}}$
has a Gaussian tail, too. For a direct proof see\ \cite{FLS}. Note that the
law of $\mathbf{B}$ is not Gaussian and there are no Fernique-type results.
\end{remark}

\begin{remark}
Potential spaces, see \cite{De04} and the references therein, are a popular
alternative to fractional Sobolev spaces. But only the latter adapt easily
to $\left( E,d\right) $-valued paths as required in rough path analysis.
\end{remark}

\begin{remark}
The $W^{\delta ,p}$-embedding of Theorem \ref{ThmfroDeltaLTone} has two
different regimes:\newline
(1) For $p$ large one has $q=1/\delta \sim 1/\alpha $. Since every $\alpha $%
-H\"{o}lder path has finite $1/\alpha $-variation (the converse not being
true) one can forget about $q$-variation.\newline
(2) When $p$ is small, the variation parameter $q=1/\delta $ can be
considerably smaller than $1/\alpha $ and $q$-variation is an essential part
of the regularity. Elementary examples show that $q$-variation does not
imply any H\"{o}lder regularity and therefore one should not forget about $%
\alpha $-H\"{o}lder regularity. The fractional Sobolev space $W^{\delta ,p}$
resp. the modulus $\omega $ are tailor-made to keep track of both regularity
aspects.
\end{remark}

\section{Cameron Martin space of fBM}

We consider fractional Brownian motion with $H\in \left( 0,1/2\right) $.
Call $\mathcal{H}^{H}$ the associated Cameron-Martin space.

\begin{theorem}
Let $1/2<\delta <$ $H+1/2$. \bigskip Then $\mathcal{H}^{H}\subset \subset
W_{0}^{\delta ,2}$.
\end{theorem}

\begin{proof}
From \cite{De04} and the references therein we know that $\mathcal{H}^{H}$
is continuously embedded in the potential space $I_{H+1/2,2}^{+}$ which we
need not define here. Then, \cite{FdeP99, De04}, $I_{H+1/2,2}^{+}\subset
W^{\delta ,2}$ so that%
\begin{equation}
\mathcal{H}^{H}\subset W^{\delta ,2}.  \label{HintoW}
\end{equation}
The compact embeddings is obtained by a standard squeezing argument: Replace 
$\delta $ by $\tilde{\delta}\in \left( \delta ,H+1/2\right) $, repeat the
argument for $\tilde{\delta}$ and then use (\ref{CompactEmbeddingWdeltap}).
\end{proof}

\begin{corollary}
For $\alpha \in \left( 0,H\right) $ and $1/\left( H+1/2\right) <q<\infty $
we have%
\begin{equation*}
\mathcal{H}^{H}\subset \subset C^{\alpha \text{-H\"{o}lder}},\,\,\,\mathcal{H%
}^{H}\subset \subset C^{q\text{-var}}.
\end{equation*}
\end{corollary}

\begin{remark}
From $\mathcal{H}^{H}\subset I_{H+1/2,2}^{+}$ it follows that $\mathcal{H}%
^{H}\subset C^{H\text{-H\"{o}lder}}$, this is well-known, \cite{De04}.
\end{remark}

\begin{remark}
For any $H\in \left( 0,1/2\right) $ we can find $1/\left( H+1/2\right) <q<2$%
. This has useful consequences. For instance, for $h,g\in \mathcal{H}^{H}$
that integral $\int hdg$ makes sense as classical Young integral with all
its continuity properties. In particular, the lift of $h\in \mathcal{H}^{H}~$%
to a geometric $p$-rough paths $p>1/H$, see \cite{MiSa05}, is well-defined
and convergence of piecewise linear approximations, uniformly over bounded
sets in $\mathcal{H}^{H}$, is an easy consequence. Such results are useful
to establish large deviations principles for enhanced Gaussian processes,
enhanced fBM being a particular example. We will discuss this in forthcoming
work.
\end{remark}

\section{Appendix}

The proof of (\ref{HintoW}) appears somewhat spread out in the references.
We present a direct argument which avoids potential spaces and fractional
calculus and extends to other Volterra kernels\footnote{%
For instance, every kernel for which one can get estimates as those in Step
2 will lead to a fractional Sobolev embedding.\newline
}.\newline
Step 1: $\mathcal{H}^{H}$ is the image of $L^{2}\left[ 0,1\right] $ under
the integral operator $K=K_{1}+K_{2}$ where 
\begin{equation*}
K_{1}\left( t,s\right) =\left( t-s\right) ^{H-1/2},\,K_{2}\left( t,s\right)
=\,s^{H-1/2}F_{1}\left( t/s\right) ;\,\,\,F_{1}=\int_{0}^{\left( \cdot
\right) -1}u^{H-3/2}\left( 1-\left( u+1\right) ^{H-1/2}\right) .
\end{equation*}%
for $s<t$ . Set $h_{i}=K_{i}g\equiv \int_{0}^{\cdot }K_{i}\left( \cdot
,s\right) g\left( s\right) ds$ with $g\in L^{2}\left[ 0,1\right] ,\,\,i=1,2.$%
\newline
Step 2: An elementary computation shows%
\begin{eqnarray*}
\sup_{u\in \left[ 0,1\right] }\int_{0}^{1-t}\left\vert K_{1}\left(
s+t,u\right) -K_{1}\left( s,u\right) \right\vert ds &=&O\left(
\,t^{H+1/2}\right) , \\
\,\sup_{s\in \left[ 0,1-t\right] }\int_{0}^{1}\left\vert K_{1}\left(
s+t,u\right) -K_{1}\left( s,u\right) \right\vert du &=&O\left(
\,t^{H+1/2}\right) .
\end{eqnarray*}%
From Cauchy-Schwartz and trivial $\sup $-estimates,%
\begin{equation*}
\left( \ast \right) =\int_{s=0}^{1-t}\left\vert h_{1}\left( s+t\right)
-h\left( s\right) \right\vert ^{2}ds=\left\vert g\right\vert
_{L^{2}}^{2}\times O\left( t^{1+2H}\right) .
\end{equation*}%
The $W^{\delta ,2}$-norm of $h_{1}$ is equivalent to $\int dt\,\left( \ast
\right) /t^{1+2\delta }$ which is less than $C\left\vert g\right\vert
_{L^{2}}^{2}$ provided $1+2H-\left( 1+2\delta \right) >-1$ and this happens
precisely for $\delta <H+1/2$.\newline
Step 3: A straight-forward computation shows (one can assume $g\in C^{1}\cap
L^{2}$ for the computation) that $\left\vert \dot{h}_{2}\right\vert
<C\left\vert g\right\vert _{L^{2}}^{2}$ provided $p<1/\left( 1-H\right) $
and hence $h_{2}\in W^{1,p}$. From (\ref{W1pEmbedding}), $W^{1,1/\left(
1-H\right) }\subset W^{H+1/2,2}$. Similarly, given $\delta <H+1/2$ we can
find $p<1/\left( 1-H\right) $, close enough to $1/\left( 1-H\right) $ so
that $W^{1,p}\subset W^{\delta ,2}$.

\newpage

\end{document}